\documentclass[12pt]{amsart}
\usepackage{amsmath,amsthm,amsfonts, amssymb}
\newtheorem{thm}{Theorem}[section]
\newtheorem{cor}[thm]{Corollary}
\newtheorem{lem}[thm]{Lemma}
\newtheorem{prop}[thm]{Proposition}

\theoremstyle{definition}

\def\e{{\rm e}}

\begin{document}

\title[Uniform ergodicity and the one-sided ergodic Hilbert transform]
{Uniform ergodicity \\and the one-sided ergodic Hilbert transform}

\author{Guy Cohen}
\address{School of Electrical Engineering, Ben-Gurion University, Beer-Sheva, Israel}
\email{guycohen@bgu.ac.il}

\author{Michael Lin}
\address{Department of Mathematics, Ben-Gurion University, Beer-Sheva, Israel}
\email{lin@math.bgu.ac.il}


\begin{abstract}
Let $T$ be a bounded linear operator on a Banach space $X$ satisfying $\|T^n\|/n \to 0$.
We prove that $T$ is uniformly ergodic if and only if the one-sided ergodic Hilbert transform
$H_Tx:= \lim_{n\to\infty} \sum_{k=1}^n k^{-1}T^k x$ converges for every $x \in \overline{(I-T)X}$.
When $T$ is power-bounded (or more generally $(C,\alpha)$ bounded for some $0< \alpha <1$),
then $T$ is uniformly ergodic if and only if the domain of $H_T$ equals $(I-T)X$.

We then study {\it rotational uniform ergodicity} -- uniform ergodicity of every $\lambda T$ 
with $|\lambda|=1$, and connect it to convergence of the rotated one-sided ergodic Hilbert 
transform, $H_{\lambda T}x$.

In the Appendix we prove that positive isometries with finite-dimensional fixed space on 
infinite-dimensional Banach lattices are never uniformly ergodic. In particular, the Koopman
operators of ergodic, even non-invertible, probability preserving transformations on  standard spaces
are never uniformly ergodic.
\end{abstract}

\dedicatory{Dedicated to Jean-Pierre Conze on his 80th birthday}
\maketitle 

\section{Introduction}
Let $T$ be a bounded linear operator on a (real or complex) Banach space $X$. We call
 $T$  {\it mean ergodic} if the averages $M_n(T):= \frac1n\sum_{k=1}^n T^k$
converge in  the strong operator topology, and {\it uniformly ergodic} if $M_n(T)$ converge
in operator norm. When $T$ is mean ergodic, the limit operator is a projection $E$ on the space 
$F(T)$ of $T$-invariant vectors, and $Ex=0$ if and only if $x \in \overline{(I-T)X}$. 
We refer to \cite[Chapter 2]{K} for additional information. Uniform ergodicity implies that 
$\sup_n\|M_n(T)\|<\infty$ ({\it Ces\`aro boundedness}) and $ \|T^n\|/n \to 0$.

Let $T$ be uniformly ergodic. Define $Y:= \overline{(I-T)X}$, which is $T$-invariant,  and put
$S:=T_{|Y}$. Then $\|M_n(S)\| \to 0$, so $\|M_n(S)\|<1$ for some large $n$ yields that $I_Y-S$ 
is invertible. Hence $Y=(I_Y-S)Y \subset (I-T)X \subset Y$, so $(I-T)X$ is closed. 

Using his operational calculus, Dunford \cite[Theorem 3.16]{D} proved that {\it $T$ on a complex
 Banach space $X$ is uniformly ergodic if and only if $\|T^n\|/n \to 0$ and $(I-T)^2X$ is closed}. 
Lin  \cite{L1} proved that {\it $T$ on a real or complex $X$ is uniformly ergodic if 
and only if $\|T^n\|/n \to 0$ and $(I-T)X$ is closed}. Mbekhta and Zem\'anek \cite{MZ} proved 
that {\it $T$ is uniformly ergodic if and only if $\|T^n\|/n \to 0$ and $(I-T)^kX$ is closed 
for some $k \ge 1$}. Kaashoek and West \cite{KW} proved that if  $\{T^n\}$ is relatively 
compact in operator norm, then $T$ is uniformly  ergodic.

When $T$ is power-bounded on $X$, we have that $T$ is uniformy ergodic if (and only if)
$\sup_n \|\sum_{k=1}^n  T^k x\| <\infty$ for every $x \in \overline{(I-T)X}$: by the 
Banach-Steinhaus theorem (for $Y$ and $S$ defined as above) we have 
$\sup_n\|\sum_{k=1}^n S^k\|<\infty$, and as before we obtain that $(I-T)X$ is closed. 
It is proved in \cite{FLR} that {\it a power-bounded $T$ is uniformly ergodic if and only if 
$\{x:\,\sup_n \|\sum_{k=1}^n T^k x\| <\infty\}$ is closed}. Examples of uniformly ergodic 
operators which are not power-bounded are given in \cite[p. 253]{LSS1}.

If $T$ satisfies $\|T^n\|/n \to 0$, then $\limsup \|T^n\|^{1/n} \le 1$, and the Abel averages
$A_r(T):=(1-r)\sum_{n=0}^\infty r^nT^n$ are well-defined (in operator norm) for $0\le r<1$. 
If $T$ is uniformly ergodic, then it is {\it Abel uniformly ergodic}, i.e. the Abel averages 
$A_r(T)$ converge in operator norm as $r\to 1^-$ \cite{H}. Conversely, if $T$ is Abel uniformly 
ergodic and $\|T^n\|/n \to 0$, then $T$ is uniformly ergodic \cite[Appendix]{L2}. 

Uniform ergodicity for $C_0$-semigroups in Banach spaces  was  studied in \cite{L2}; 
the study was extended to some Fr\'echet spaces by Albanese et al. \cite{ABR}.

\section{Uniform ergodicity and the one-sided ergodic Hilbert transform}

Let $T$ be a bounded linear operator on $X$, and put  $H_nx:=\sum_{k=1}^n \frac1k T^k x$.
Let $D(H)$  be the set of $x$ for which the limit $H_Tx:=\lim_{n\to \infty} H_nx$ exists; the 
operator $H=H_T$ is called {\it the one-sided ergodic Hilbert transform} (of $T$). If $x \in D(H)$, 
then for every $T^*$-invariant functional $x^* \in X^*$ we have $\langle x^*,x\rangle =0$, so
 by the Hahn-Banach theorem $x \in \overline{(I-T)X}$; hence $D(H) \subset \overline{(I-T)X}$.
The following was proved by Becker \cite{B2}.

\begin{lem} \label{inclusion}
 If $\sup_n \|M_n(T)\| < \infty$ and $\frac{T^n }n \to 0$ strongly,
then ${(I-T)X \subset D(H) }$. 
\end{lem}
Becker \cite{B2} has an example of $T$ with $n^{-1}\|T^n\| \to 0$ such that $(I-T)X \not\subset D(H)$.
 \begin{cor} \label{inclusions}
If $T$ is power-bounded or mean ergodic, then $(I-T)X \subset D(H) \subset\overline{(I-T)X}$.
\end{cor}

\begin{prop} \label{coboundaries}
If $T$ is uniformly ergodic, then $(I-T)X=D(H)=\overline{(I-T)X}$.
\end{prop}
\begin{proof} 
Since $T$ is uniformly ergodic,  $(I-T)X$ is closed, so Corollary \ref{inclusions} 
yields $(I-T)X=D(H)=\overline{(I-T)X}$.
\end{proof}
An application of Proposition \ref{coboundaries} yields a new proof of the following well-known result 
(e.g. \cite[Theorem 13]{CL} and the discussion following it).

\begin{prop} \label{halmos}
Let $\theta$ be an ergodic invertible measure preserving transformation on an atomless probability
space $(\mathcal S, \Sigma,\mu)$, and let $Tf:=f\circ \theta$ be the induced unitary operator on 
$L_2(\mu)$. Then $T$ is not uniformly ergodic.
\end{prop}
\begin{proof} Halmos \cite{Ha} observed that $(I-T)L_2(\mu) \subset D(H)$, and proved that 
there exists $f $ in $\overline{(I-T)L_2(\mu)}$ for which the series 
	$\sum_{n=1}^\infty \frac1n T^nf$ 
does not converge in $L_2$, i.e. $f \notin D(H)$. Hence $(I-T)L_2(\mu)$ is not closed, and 
$T$ is not uniformly ergodic.
\end{proof}

{\bf Remarks.} 1. Proposition \ref{halmos} implies that  $T$ induced by the above $\theta$,
has no uniform rate of convergence in von Neumann's mean ergodic theorem (see
\cite[Theorem 13]{CL}, or the proof of \cite[Proposition 3.2]{DL2}).

2. The operator $Tf=f\circ \theta$ of Proposition \ref{halmos} is an invertible isometry of each
$L_p(\mu)$, $1 \le p< \infty$. Dowker and Erd\"os \cite[Example 2]{DE} proved the existence of 
$f \in L_\infty$  with $\int f\,d\mu =0$ such that $H_nf$ does not converge in $L_p$
for any $1\le p<\infty$, so $T$ is not uniformly ergodic in any $L_p(\mu)$;
moreover the series $\sum_{n=1}^\infty \frac1n T^nf$ almost everywhere does not converge.
See also \cite[Corollary 3.5]{DL2}.

3. Without ergodicity, Proposition \ref{halmos} may fail; e.g. $\theta x=x$.

4. For a proof of Proposition \ref{halmos} without assuming invertibility, see Corollary 
\ref{no-inverse} below.

5. Derrienic and Lin \cite[Proposition 2.20]{DL} proved that when $T$ is power-bounded (and not 
uniformly ergodic), for $0<\alpha <1$ we have $(I-T)^\alpha X\subset D(H)$; if in addition $T$ 
is mean ergodic. then $x \in (I-T)^\alpha X$ satisfies $\|M_n(T)x\| =o(1/n^{1-\alpha})$.
Gomilko, Haase and Tomilov \cite[Theorem 4.1]{GHT} proved that if $T$ is power-bounded and 
$x \in D(H)$, then $\|M_n(T)x\|=o(1/\log n)$; of course, this is of interest only when $T$ is not 
uniformly ergodic, in which case $D(H) \ne (I-T)X$ (see Theorem \ref{DL} below), and 
$D(H) \ne \overline{(I-T)X}$ by \cite[Theorem 2.23]{DL} (see also Theorem \ref{becker}(ii) below).

6. Let $\{\lambda_k\}_{k \ge 1} \subset \mathbb T$ satisfy $\delta := \inf_k |1- \lambda_k|>0$. 
On a separable infinite-dimensional complex Hilbert with orthonormal basis $\{e_k\}_{k\ge 1}$
define $T(\sum_k a_k e_k)= \sum_k \lambda_k a_k e_k$. Then $T$ is unitary with $I-T$ invertible,
so $T$ is uniformly ergodic. The same construction yields a uniformly ergodic invertible isometry on 
the complex $\ell_p(\mathbb N)$, $1\le p<\infty$, when $e_k:=1_{\{k\}}$ is the $k$th unit vector.

7. Necessary and sufficient conditions for the convergence of $H_nx$ when $T$ is an isometry or a
normal contraction in a Hilbert space are given in \cite[Theorem 2.2.1]{CCC}.
\medskip

Cohen, Cuny and Lin \cite[Theorem 3.9]{CCL} proved that if $T$ is power-bounded on $X$, then $-H$ 
is the infinitesimal generator of the semigroup $\{(I-T)^r: r\ge 0\}$ defined on $\overline{(I-T)X}$
(see \cite[Theorem 2.22]{DL}). Hence $H$ is a closed operator, with $D(H)$ dense in $\overline{(I-T)X}$.
Haase and Tomilov \cite[Theorem 6.2]{HT} proved the same result independeltly, and identified $H$ 
with $-\log(I-T)$. This identification was extended to $(C,\alpha)$-bounded operators, $0<\alpha<1$,
by Abadias, Gal\'e and Lizama \cite[Theorem 8.10]{AGL}.

\begin{thm} \label{becker}
Let $T$ be a bounded linear operator on a (real or complex) Banach space $X$. Then the following are
equivalent:

(i) $T$ is uniformly ergodic.

(ii) $\frac1n\|T^n\| \to 0$ and $D(H)= \overline{(I-T)X}$.

(iii)  $\frac1n\|T^n\| \to 0$ and $\sup_n \|H_n x\| < \infty$ for every $x \in \overline{(I-T)X}$.

(iv) $\sup_n \|M_n(T)\| < \infty$, $\ \frac1n\|T^n\| \to 0$, and $D(H)$ is closed.

\end{thm}
\begin{proof} It is shown in Proposition \ref{coboundaries} that (i) implies (ii) and (iv).

(iv) implies (ii): By \cite[Propositions 2.1 and 2.2]{B2}, 
$(I-T)X \subset D(H) \subset \overline{(I-T)X}$, so $D(H)= \overline{(I-T)X}$.

Clearly (ii) implies (iii).

(iii) implies (i): We adapt the proof of uniform ergodicity in \cite[Theorem]{B1}.
 Let $\tilde H_n$ be the restriction of $H_n$ to $Y=\overline{(I-T)X}$. 
By the Banach-Steinhaus theorem, $\sup_n \|  \tilde H_n\| = K < \infty$. 
By Abel's summation formula (or by an easy induction), $S:=T_{|Y}$ satisfies
$$
M_n(S)+\frac{n+1}n \tilde H_n -\frac1n \sum_{k=1}^n \tilde H_k,
$$
which yields that $\sup_n \|M_n(S)\| = C < \infty$. For fixed $j>0$ we have
$$
S^j\sum_{k=1}^n S^k - \sum_{k=1}^n S^k=\sum_{k=n+1}^{n+j} S^k -\sum_{k=1}^j S^k=
S^n jM_j(S) - jM_j(S).
$$
Since $n^{-1}\|S^n\| \to 0$, we obtain $\lim_{n\to \infty}\|S^jM_n(S)-M_n(S)\| =0$.

Define $B_n:=\tilde H_n/\sum_{k=1}^n \frac1k$. Then $\|B_n\| \to 0$.
Fix $\varepsilon >0$, and choose $k$ large so $\|B_k\|<\varepsilon$. Since
$B_k$ is a convex combination of $S, \dots,S^k$, also $\lim_n\|B_kM_n(S)=M_n(S)\|=0$.
Then for $n>N$ we have $\|B_kM_n(S)=M_n(S)\|< \varepsilon$, and
$$
\|M_n(S)\| \le \|M_n(S)-B_kM_n(S)\|+\|B_kM_n(S)\| <\varepsilon(1+C).
$$
Taking $\varepsilon <1/(1+C)$ we obtain $\|M_n(S)\| <1$ for large $n$, so $I-M_n(S)$ is invertible 
on $Y$, hence also $I-S$ is inverible on $Y$, which implies $(I-S)Y=Y$. Then
$Y=(I-S)Y \subset (I-T)X \subset Y$, which shows  that $(I-T)X$ is closed, and
by \cite{L1} $T$ is uniformly ergodic.
\end{proof}

{\bf Remarks.} 1. If in (ii) or (iv) of Theorem \ref{becker} we replace $\frac1n \|T^n\| \to 0$
by $\frac1n T^n \to 0$ strongly, then $T$ need not be uniformly ergodic, even when it is
mean ergodic. In Burlando's \cite[Example 3.8]{Bu}, $T$ is mean ergodic with $(I-T)X$ closed,
but $T$ is  not uniformly ergodic; by Corollary \ref{inclusions} 
$(I-T)X \subset D(H) \subset \overline{(I-T)X}$, so $D(H)$ is closed and equals 
$(I-T)X=\overline{(I-T)X}$.

2. The above proof of (iii) implies (i) yields that if $T$ is uniformly ergodic, then $H$ is a bounded 
operator on $(I-T)X=\overline{(I-T)X}$.

\begin{cor} \label{laura}
	$\|\frac1n\sum_{k=1}^{n} T^k\| \to 0$ if and only if $\|T^n\|/n \to 0$ and $D(H)=X$.
\end{cor}
\begin{proof}  When $\|M_n\| \to 0$ we have $I-T$ invertible, and then $D(H)=(I-T)X=X$ using 
	Proposition \ref{coboundaries}.

When $D(H)=X$ we have $D(H)=\overline{(I-T)X}$, since $D(H) \subset \overline{(I-T)X}$, 
so $T$ is uniformly ergodic by Theorem \ref{becker}(ii). Then $(I-T)X=X$, and $\|M_n\| \to 0$.
\end{proof}

{\bf Remarks.} 1. The assumption $D(H)=X$ alone always implies, by Kronecker's lemma for vectors,  
that $T$ is mean ergodic, with $M_nx \to 0$ for every $x \in X$.

2. Let $T$ be the mean ergodic operator in Burlando's example mentioned above, and put
 $X_0=\overline{(I-T)X}=(I-T)X$, $T_0=T_{|X_0}$. Then $D(H_{T_0})=X_0$ (since $H_T=H_{T_0}$),
 but $T_0$ is not uniformly ergodic. Thus the assumption $\|T^n\|/n \to 0$ in Corollary \ref{laura}
 cannot be weakened to $T^n/n \to 0$ strongly.


\begin{cor} \label{series}
If $\sum_{n=1}^N  n^{-1} T^n(I-T)$ converges in operator norm as $N \to \infty$, and $D(H)$ is
 closed, then $T$ is uniformly ergodic.
\end{cor}
\begin{proof}
By Kronecker's lemma for vectors, the  convergence of the series  implies \newline
${\|\frac1n\sum_{k=1}^n T^k(I-T)\| \to 0}$. Hence
$$
\Big\|\frac{I-T^n}n \Big\| = \Big\|\frac1n\sum_{k=0}^{n-1} T^k(I-T)\Big\| \le
\frac{\|I-T\|}n+\frac{n-1}n \Big\|\frac1{n-1}\sum_{k=1}^{n-1} T^k(I-T)\Big\| \to 0,
$$
so $\|T^n\|/n \to 0$. Since $H_n\cdot (I-T)=\sum_{k=1}^n k^{-1}{T^k(I-T)}$, we obtain
$(I-T)X \subset D(H)$. Since always $D(H)  \subset \overline{(I-T)X}$, we coclude that 
$D(H)=\overline{(I-T)X}$ when $D(H)$ is closed. Now Theorem \ref{becker}(ii) yields uniform ergodicity.
\end{proof}

\begin{cor} \label{kt}
If $\|T^n(I-T)\|=O(1/(\log n)^\delta)$ for some $\delta>1$ and $D(H)$ is closed, then $T^n$
converges in operator norm.
\end{cor}
\begin{proof} 
$T$ is uniformly ergodic by Corollary \ref{series} and $\|T^n(I-T)\| \to 0$, so $T^n$ converges 
in norm by \cite[Theorem 6]{Lu}.
\end{proof}

\begin{cor} \label{log}
If $\|T^n\|=O(n/(\log n)^\delta)$ for some $\delta>1$ and $D(H)$ is closed, then $T$ is uniformly
ergodic.
\end{cor}
\begin{proof}
Clearly $\|T^n\|/n \to 0$. Application of the growth assumption to the last series below
$$\sum_{k=1}^{n} \frac{T^k(I-T)}k=  T -\frac{T^{n+1}}n - \sum_{k=2}^n \frac{T^k}{k(k-1)}\ $$
yields that the series on the left converges in norm, and then Corollary \ref{series} applies.
\end{proof}

{\bf Remarks.} 1. The proof of Corollary \ref{series} shows that 
$(I-T)X \subset D(H) \subset \overline{(I-T)X}$ when  $\|T^n\|=O(n/(\log n)^\delta)$ with $\delta>1$.
Becker \cite{B2} has an example of $T$ with $n^{-1}\|T^n\| \to 0$ such that $(I-T)X \not\subset D(H)$
 (so $T$ is not even mean ergodic, by Corollary \ref{inclusions}); in fact, $\sup_n \|M_n(T)\|=\infty$.
 When in her example $a_j=(j+1)/\log(j+1)$, we have $\|T^n\| \sim O(n/\log n)$. 

2. Corollary \ref{log} implies the result of \cite{B1}: {\it  If $n^{-\alpha}\|T^n\| \to 0$ for some 
$\alpha \in [0,1)$ and $D(H)$ is closed, then $T$ is uniformly ergodic.} 

\medskip

 Proposition  \ref{coboundaries} raises the question whether a (Ces\`aro bounded) operator $T$, 
satisfying $\ \frac1n\|T^n\| \to 0$ and $D(H)=(I-T)X$, must be uniformly ergodic.

\begin{lem} \label{growth-eht}
Let $y$ satisfy $\|M_n(T)y\|=O(1/n^\beta)$ for some $0< \beta <1$. Then $y \in D(H)$.
\end{lem}
\begin{proof}
Put $S_k:=\sum_{j=1}^k T^j$ and $\alpha=1-\beta$. Then $\|S_ny\|=O(n^\alpha)$.
For $n>\ell>1$ we have
$$
\sum_{k=\ell}^n \frac{T^k y}k=\sum_{k=\ell}^n\frac{S_{_k}y-S_{_{k-1}}y}k =
\frac{S_{_n}y}n +\sum_{k=\ell}^{n-1}\Big(\frac1k-\frac1{k+1}\Big)S_ky-\frac{S_{_{\ell-1}}y}\ell.
$$
Since $S_n y/n \to 0$ and
$$
\Big\|\sum_{k=\ell}^{n-1} \Big(\frac1k -\frac1{k+1} \Big) S_{_k}y \Big\| =
\Big\|\sum_{k=\ell}^{n-1} \Big(\frac1{k(k+1)} \Big) S_{_k}y \Big\| \le
C \sum_{k=\ell}^{\infty} \Big(\frac1{k^{1-\alpha}(k+1)} \Big),
$$
the sequence $\big\{ \sum_{k=1}^n \frac{T^k y}k \big\}$ is Cauchy, so $y \in D(H)$.
\end{proof}

\begin{thm} \label{DL}
Let $T$ be power-bounded on a (real or complex) Banach space $X$.  
Then the following are equivalent:

(i) $T$ is uniformly ergodic.

(ii) $D(H)$ is closed.

(iii) $D(H)=\overline{(I-T)X}$.

(iv) $D(H)=(I-T)X$.

(v) For every $y \in \overline{(I-T)X}$ there is $0<\beta<1$ such that $\|M_n(T)y\| =O(1/n^\beta)$.
\end{thm}
\begin{proof}
The equivalence of the first three conditions follows from Theorem \ref{becker}.
By Proposition \ref{coboundaries}, (i) implies (iv).

Let $T$ be power-bounded with $D(H)=(I-T)X$. Assume $T$ is not uniformly ergodic, and
 fix $\alpha \in(0,1)$. Derriennic and Lin \cite{DL} used the Taylor series of 
$(1-t)^\alpha=1 -\sum_{j=1}^\infty a_j^{(\alpha)}t^j$ to define 
$(I-T)^\alpha: =I -\sum_{j=1}^\infty a_j^{(\alpha)} T^j$, with convergence in operator norm.
It is proved in \cite[Propostion 2.2]{DL} that if $(I-T)X$ is not closed, then 
$(I-T)X \subset (I-T)^\alpha X$ with no equality. By \cite[Proposition 2.20]{DL}, 
$(I-T)^\alpha X \subset D(H)$, which contradicts $D(H)=(I-T)X$; hence $(I-T)X$ is closed,
and $T$ is uniformly ergodic by \cite{L1} .

When (i) holds, $\overline{(I-T)X}=(I-T)X$, so $\|M_n(T)y\|=O(1/n)$ for $y \in\overline{(I-T)X}$.

	Assume (v). If $y \in \overline{(I-T)X}$, then $y \in D(H)$ by Lemma \ref{growth-eht}.
Thus $\overline{(I-T)X} \subset D(H)$, hence $D(H)= \overline{(I-T)X}$, and by Theorem 
\ref{becker}(ii) $T$ is uniformly ergodic.
\end{proof}

{\bf Remark.}  When $T$ is unitary on a complex Hilbert space, Kachurovskii et al. \cite[Theorem 1]{KPK}
gave a spectral characterization of $\|M_n(T)y\|= O(1/n^\beta)$ when $y \in \overline{(I-T)X}$ and
$\beta <1$.

\begin{lem} \label{exponential}
	Assume that $\|T^n\| \to 0$. Then there exist $C>0$ and $r<1$ such that $\|T^n\| \le Cr^n$.
	Hence $I-T$ is invertible with $(I-T)^{-1}=\sum_{n=0}^\infty T^n$ (convergence in operator norm).
\end{lem}
The first part is proved (for $X$ real or complex) in \cite{DL2}, and yields the second part.
\smallskip

{\bf Remark.} When $T$ is a positive operator on a Banach lattice $X$ and $\|T^n\| \to 0$, 
the lemma yields that $(I-T)^{-1}$ is positive. If $X$ is a complex  Banach lattice, $T$ is positive 
with $r(T)\le 1$, and $(I-T)$ is invertible, then $1 \notin \sigma(T)$, and since $r(T) \in \sigma(T)$
\cite[p. 323]{Sh}, we have $r(T)<1$; hence $\|T^n\| \to 0$. 
Gl\"uck and Mironchenko  \cite[Theorem 3.3]{GM}
 proved that when $X$ is a complex  Banach lattice  and $T$ is positive, $\|T^n\| \to 0$ if $I-T$
is invertible and its  inverse ia positive (without any assumption on $r(T)$).

\begin{lem} \label{invert}
Let $T$ be (weakly) mean ergodic. If $x \in (I-T)X$, then 
$\frac1N\sum_{n=1}^N\sum_{k=0}^{n-1} T^k x $
	converges (weakly), say to $y$, $\ y \in \overline{(I-T)X}$, and $(I-T)y=x$.
\end{lem}
\begin{proof}
By the ergodic decomposition, there exists a unique $y \in \overline{(I-T)X}$ with $(I-T)y=x$.
	Then $M_ny \to 0$ (weakly), and
$$
 \frac1N\sum_{n=1}^N\sum_{k=0}^{n-1} T^k x =
 \frac1N\sum_{n=1}^N\sum_{k=0}^{n-1} T^k (I-T)y=\frac1N\sum_{n=1}^N(I-T^n)y=
	y-M_Ny \underset{N\to\infty}\to y \quad \text{(weakly)}.
$$
\end{proof}


{\bf Remark.} If $T$ is (weakly) mean ergodic and $I-T$ is invertible, Lemma \ref{invert} yields
that $\frac1N\sum_{n=1}^N\sum_{k=0}^{n-1} T^k \to (I-T)^{-1} $ (weakly).
If $T$ is mean ergodic, then $D(H)=X$ by Corollary \ref{inclusions}, but $T$ need not be 
uniformly ergodic; see remarks following Corollary \ref{laura}.

\begin{cor}
	Assume that $\|M_n\| \to 0$. Then $I-T$ is invertible, and
$\frac1N\sum_{n=1}^N\sum_{k=0}^{n-1} T^k \to (I-T)^{-1}$ in operator norm.
\end{cor}
\noindent {\it Proof.} The assumption implies that $I-T$ is invertible (see the introduction). 
Then, by the above algebra,
$$
 \frac1N\sum_{n=1}^N\sum_{k=0}^{n-1} T^k  =
(I-T)^{-1}\frac1N\sum_{n=1}^N(I-T^n)=(I-T)^{-1}(I-M_N) \to (I-T)^{-1}. \qquad \qquad  \square
$$


\begin{prop} \label{kalton}
The following are equivalent for $T$ on a complex Banach space $X$:

(i) $T^n$ converges in operator norm.

(ii)  $T$ is uniformly ergodic and $\|T^n(I-T)\| \to 0$.

(iii) $\limsup_{n\to\infty} n \|T^n(I-T) \| < \e^{-1}$.

(iv) $\|T^n(I-T)\| \le Cr^n$ for some $C>0$ and $r<1$.

	(v) $r(T_{|\overline{(I-T)X}})<1$.
\end{prop}
\begin{proof}
The equivalence of (i) and (ii) is a corollary of \cite[Theorem 6]{Lu}.


By the ergodic decomposition (i) implies (v). Assume (v) and put $Y=\overline{(I-T)X}$ and $S=T_{|Y}$.
	Then $r(S)<1$ inplies $\|S^n\| \to 0$, and 
	$\|T^n(I-T)\| \le\|S^n\|\cdot\|I-T\| \le \|I-T\|C_0r^n$ for some $r<1$ 
	by Lemma \ref{exponential}, so (iv) holds.

Clearly (iv) implies (iii);  Kalton et al. \cite[Theorem 3.1]{KMOT} proved that (iii) implies (i).
\end{proof}
{\bf Remarks.} 1. In conditions (ii) -- (v) there is no a priori assumption of power-boundedness, 
which follows from (i).

2. If $\limsup_n n \|T^n(I-T) \| \ge \e^{-1}$, then $T$ need not be power-bounded
\cite[Theorem 3.6]{KMOT}.

3. See \cite[Theorem 5]{Ja} for a different proof of (v) implies (i).

\begin{prop} 
Let $T$ be power-bounded on $X$ complex. Then the following are equivalent:

(i) $T^n$ converges in operator norm.

(ii) $\|T^n(I-T)\| \to 0$ and $D(H)$ is closed.

(iii)  $\sigma(T)\cap\mathbb T \subset\{1\}$ and $D(H)$ is closed. 

\end{prop}
\begin{proof}
 By Theorem \ref{becker}, (i) implies (ii).

By Katznelson and Tzafriri \cite{KT}, $\sigma(T)\cap\mathbb T \subset\{1\}$ is  equivalent to  
$\|T^n(I-T)\| \to 0$ (for $T$ power-bounded). Hence (ii) and (iii) are equivalent.

Assume (ii). Since $D(H)$ is closed, Corollary \ref{inclusions} and power-boundedness yield
 $D(H)=\overline{(I-T)X}$, so by Theorem  \ref{becker}(ii) $T$ is uniformly ergodic.  
Norm convergence  of $T^n$ follows from \cite{Lu}.
\end{proof}

{\bf Remarks.} 1. Power-boundedness is necessary for convergence of $T^n$.

2. Conditions (i) and (ii) are equivalent also for $T$ power-bounded on a real Banach space.

3. In conditions (ii) and (iii) of the theorem,  $D(H)$ closed can be replaced by $D(H)=(I-T)X$,
and uniform ergodicity is then obtained from Theorem \ref{DL}.

4. $T$ power-bounded with $\sigma(T)\cap \mathbb T=\{1\}$ need not be uniformly ergodic, 
though ${\|T^n(I-T)\| \to 0}$. See Example following Proposition \ref{eht}.

 5. By \cite[Lemma 2.1]{LSS1}, if $T$ is uniformly ergodic with $\sigma(T)=\{1\}$, then $T=I$.

\begin{thm} \label{gluck}
Let $T$ be a bounded linear operator on a complex Banach space. Then ${\|T^n\| \to 0}$ if 
(and only if) for every $x \in X$ and $x^* \in X^*$ there is a $p \in [1,\infty)$ such that 
\begin{equation} \label{ell-p}
\sum_{n=1}^\infty |\langle x^*,T^n x\rangle|^p < \infty. 
\end{equation}
\end{thm}
\begin{proof}
If $\|T^n\|  \to 0$, then the convegence is exponentially fast, so \eqref{ell-p} holds for 
every $p\ge 1$.

Assume that \eqref{ell-p} holds for every $x \in X$ and $x^* \in X^*$ ($p$ may depend on $x$ and
 $x^*$). Then $T^n \to 0$ in the weak operator topology, so $T$ is power-bounded, and has
spectral radius $r(T) \le 1$. 
 When \eqref{ell-p} holds,  we may assume $p>1$. By H\"older's inequality ($q=p/(p-1)$),
$$
\sum_{n=1}^\infty \frac{|\langle x^*,T^n x\rangle|}n \le 
\Big(\sum_{n=1}^\infty \frac1{n^q}\Big)^{1/q}
\cdot(\sum_{n=1}^\infty |\langle x^*,T^n x\rangle|^p)^{1/p} < \infty.
$$
Applying  the Banach-Steinhaus theorem, we obtain that $\sup_n \|H_nx\|< \infty$ for every
$x \in X$, so by Theorem \ref{becker}(iii) $T$ is uniformly ergodic, with $\|M_n(T)\| \to 0$. Hence
$I-T$ is invertible. Let $\lambda \in \mathbb T$.  Then also $\lambda T$ satisfies \eqref{ell-p},
so $\bar \lambda I-T$ is invertible. Since  $r(T) \le 1$, and by the above  
$\sigma(T) \cap \mathbb T=\emptyset$, we obtain $r(T)<1$, so $\|T^n\| \to 0$.
\end{proof}
{\bf Remarks.} 1. Theorem \ref{gluck} was first proved, differently,  by Gl\"uck \cite{Gl}

2. For $p$ constant, independent of $x$ and $x^*$, Theorem \ref{gluck} was proved by G. Weiss
 \cite{We}.

\begin{cor} \label{weiss}
Let $T$ on $X$ complex satisfy $\|T^n\|/n \to 0$. If for every $x \in \overline{(I-T)X}$ 
there is a $p \in (1,\infty)$ such that
$\sum_{n=1}^\infty \|T^n x\|^p < \infty$, then $T^n$ converges in operator norm.
\end{cor}
 \begin{proof}
Put $Y=\overline{(I-T)X}$ and $S=T_{|Y}$. Then $S$ satisfies \eqref{ell-p} for every $x \in Y$ 
and $x^* \in Y^*$, so by Theorem \ref{gluck} $\|S^n\| \to 0$. By the previous proof 
  $D(H)=\overline{(I-T)X}$, so by Theorem \ref{becker}  $T$ is uniformly ergodic, with the limit $E$
projecting on $F(T)$ according to the ergodic decomposition $X=F(T) \oplus Y$. Then
$\|T^n-E \| =\|T^n(I-E)\| \le \|S^n\|\cdot\|I-E\| \to 0$.
\end{proof}
{\bf Remark.}  When $X$ is a Hilbert space and $p$ is fixed, Corollary \ref{weiss} yields the
discrete version of the Datko-Pazy theorem for $C_0$-semigroups.
\bigskip

\section{$(C,\alpha)$ uniform ergodicity}

In this section $T$ is a bounded linear operator on a Banach space $X$. When $X$ is over $\mathbb C$,
we denote by $R(\lambda,T)$ the resolvent of $T$.

For $\beta>-1$ put $A_0^\beta=1$ and for $n\ge 1$  put 
$A_n^\beta =(\beta+n)(\beta +n-1)\cdots(\beta+1)/n!$
\smallskip

{\bf Definition.} The {\it Ces\`aro means of order $\alpha$} of an operator $T$ are defined by
$M_n^\alpha =S_n^\alpha/A_n^\alpha$, where 
$S_n^{\alpha}=\sum_{k=0}^n A_{n-k}^{\alpha-1} T^k$. We have $M_n^0=S_n^0=T^n$
and $M_n^1=\frac1{n+1}\sum_{k=0}^n T^k$.

We say that $T$ is $(C,\alpha)$ {\it uniformly ergodic} if $M_n^\alpha$ converges in operator norm,
say to $E$.

As shown below (Proposition \ref{summability} or Corollary \ref{ue}), $(C,\alpha)$ uniform ergodicity
with ${0<\alpha<1}$ implies uniform ergodicity.
\smallskip

When $X$ is a complex space, Hille \cite[Theorem 6]{H} proved that necessary conditions for 
$(C,\alpha)$ uniform ergodicity, $\alpha>0$, are $\|T^n\|/n^\alpha \to 0$ and 
$\|(\lambda-1)R(\lambda,T)-E\| \underset{\lambda \to 1^+} \to 0$.
Ed-dari \cite{Ed} proved that these two conditions together imply $(C,\alpha)$ uniform ergodicity.
Yoshimoto \cite{Y} proved that in that case $(I-T)X$ is closed, $X=F(T)\oplus (I-T)X$, and $E$
projects on $F(T)$ according to this decomposition.
Badiozzaman and Thorpe \cite[Theorem 2]{BT} and (independently) Yoshimoto \cite[Theorems 1,3]{Y} 
proved that {\it  $T$ is $(C,\alpha)$ uniformly ergodic if and only if it is Abel uniformly ergodic 
and $\|T^n\|/n^\alpha \to 0$}. 
When $T$ is power-bounded, uniform ergodicity implies $(C,\alpha)$ uniform ergodicity for every 
$\alpha>0$ \cite[Theorem 3]{H}.  Combining Hille's and Ed-dari's results, we obtain the following.

\begin{prop} \label{summability}
Fix $\alpha>0$. Let $T$ on a complex Banach space be $(C,\alpha)$ uniformly ergodic. Then $T$
 is $(C,\beta)$ uniformly ergodic for every $\beta > \alpha$. In particular, when $\alpha<1$ $\ T$ 
is uniformly ergodic.
\end{prop}

The following proposition is probably well-known, but we could not find a direct reference 
(it is proved in \cite{H} when $M_n^\alpha$ converges uniformly).
\begin{prop} \label{growth}
Let $0<\alpha<1$ and let $T$ on a (real or complex) Banach space $X$ be $(C,\alpha)$-bounded, 
i.e. $\sup_n \|M_n^\alpha\| =M < \infty$. Then $\|T^n\| =O(n^\alpha)$.
\end{prop}
\begin{proof}
By \cite[(1.10), p. 77]{Z}, $T^n=S_n^0=\sum_{k=0}^n A_{n-k}^{-\alpha-1} S_k^\alpha$. Then
$$
\frac{S_n^0}{A_n^\alpha} = 
\frac{1}{A_n^\alpha}\sum_{k=0}^n A_{n-k}^{-\alpha-1} A_k^{\alpha} M_k^{\alpha}.
$$
For $\alpha>0$ we have $A_n^\alpha >0$ and increasing \cite[Theorem 1.17]{Z}, so
$$
\Big\|  \frac{S_n^0}{A_n^\alpha} \Big \| \le M \sum_{k=0}^\infty |A_{k}^{-\alpha-1}|. 
$$
By \cite[(1.16), p. 77]{Z} the infinite series converges, and by definition, 
$A_n^\alpha \sim n^\alpha$.  We conclude that $\|T^n\|/n^\alpha$ is bounded.
\end{proof}

\begin{thm} \label{C-alpha}
Fix $0<\alpha<1$. An operator  $T$ on a (real or complex) Banach space is $(C,\alpha)$ uniformly
 ergodic if and only if $(I-T)X$ is closed and $\|T^n\|/n^\alpha \to 0$.
\end{thm}
\begin{proof} 
 If $(I-T)X$ is closed  and $\|T^n\|/n^\alpha \to 0$, then $T$ is uniformly ergodic by \cite{L1}
 with limit $E$, and $X=F(T) \oplus (I-T)X$; the limit operator $E$ projects onto $F(T)$ with 
$\ker E=(I-T)X$.
By Lemma 1 of \cite{Y}, which is valid also in the real case, $\lim_n \|M_n^\alpha(I-T)\| = 0$.
Since $(I-T)$ is invertible on $(I-T)X$ and $M_n^\alpha x=x$ for $x \in F(T)$, we obtain
$$
\|M_n^\alpha-E\| = \|M_n^\alpha(I-E)\| \le 
$$
$$
\|I-E\| \cdot\|M_n^\alpha(I-T)\| \cdot\|(I-T_{|(I-T)X})^{-1}\| \to 0.
$$
\smallskip

Assume now that $M_n^\alpha$ is uniformly ergodic, with limit $E$. As indicated in
 \cite[Theorem 1]{H}, the proof that $\|T^n\|=o(n^\alpha)$ follows verbatim from the proof for
 numerical sequences. The proof of the general \cite[Theorem 1]{H} does not use the assumption 
that $X$ is over $\mathbb C$, and we coclude that $T$ is uniformly Abel ergodic (with the same 
limit $E$). Since $\alpha<1$, by \cite[Appendix]{L2}, $(I-T)X$ is closed. 
\end{proof}

{\bf Remarks.} 1. The case $\alpha=1$ is proved in \cite{L1}  (and used in the above proof).

 2. For $X$ complex, the theorem follows from combining \cite{BT} with  \cite[Appendix]{L2}.

\begin{cor} \label{ue}
Fix $0<\alpha <1$. If $T$ on $X$ real or complex  is $(C,\alpha)$ uniformly ergodic, then it is 
$(C,\beta)$ uniformly ergodic for any $\alpha<\beta \le 1$.
\end{cor}

\begin{lem} \label{me}
Let $0<\alpha<1$ and let $T$ be $(C,\alpha)$-bounded on a  Banach space $X$. 
Then $\sup_n \|M_n(T)\|< \infty$ and $\frac1n \|T^n\| \to 0$.
\end{lem}
\begin{proof} Since $0<\alpha<1$, the lemma follows from \cite[Lemma 1]{De} with $\beta=1$.
\end{proof}

\begin{cor}
Let $T$ be $(C,\alpha)$ bounded on $X$, $0<\alpha <1$. Then $(I-T)X \subset D(H)$.
\end{cor}
\begin{proof} Combine Lemma \ref{me} with Lemma \ref{inclusion}.
\end{proof}

\begin{cor} \label{beta-ge-alpha}
Let $T$ be $(C,\alpha)$ bounded on $X$, $0<\alpha <1$. If $D(H)$ is closed, then $T$ is uniformly
ergodic; moreover, $T$ is $(C,\beta)$ uniformly ergodic for any $\beta >\alpha$.
\end{cor}
\begin{proof} By Lemma \ref{me} and Theorem \ref{becker}(iv), or by Proposition \ref{growth}
and Corollary \ref{log}, $T$ is uniformly ergodic.
Hence $(I-T)X$ is closed, and $X=F(T) \oplus (I-T)X$. By \cite[Lemma 1]{De}, for $\beta>\alpha$
the $(C,\beta)$ averages $M_n^\beta$ satisfy $\|(I-T) M_n^\beta \| \to 0$. Together with 
$ M_n^\beta x=x$  for $x \in F(T)$, we obtain that $M_n^\beta$ converges uniformly.
\end{proof}

{\bf Remark.} Li, Sato and Shaw \cite[Proposition 4.4]{LSS} showed that on any 
infinite-dimensional $X$ there exists $T$ which is Ces\`aro bounded, but not $(C,\alpha)$ 
bounded for any $\alpha \in [0,1)$.

\begin{cor}
Fix $0<\alpha<1$.
 If $\|T^n\|/n^\alpha \to 0$ and $D(H)$ is closed, then $T$ is $(C,\alpha)$ uniformly ergodic.
\end{cor}
\begin{proof}
By Corollary \ref{log}, $T$ is uniformly ergodic, so $(I-T)X$ is closed.
Hence $T$ is $(C,\alpha)$ uniformly ergodic by  Theorem \ref{C-alpha}.
\end{proof}

The next proposition extends Theorem \ref{DL}(v) to $(C,\alpha)$ uniform ergodicity.
\begin{prop}
Fix $0< \alpha <1$. The operator $T$ is $(C,\alpha)$ uniformly ergodic if and only if 
$\|T^n\|/n^\alpha \to 0$ and for every $y \in \overline{(I-T)X}$ there is $0<\beta<1$ such that
 $\|M_n(T)y\| =O(1/n^\beta)$.
\end{prop}
\begin{proof}
Let $T$ be $(C,\alpha)$ uniformly ergodic. Then by Theorem \ref{C-alpha}, $\|T^n\|/n^\alpha \to 0$
and $(I-T)X$ is closed. 
For $y \in \overline{(I-T)X}= (I-T)X$ there is $x$ with $y=(I-T)x$. Since $\|T^n\|/n^\alpha$ is bounded,
 we have $\|\sum_{k=0}^{n-1}T^ky\| \le \|x\|(1 +\|T^n\|) = O(n^\alpha)$, so 
$\|M_n(T)y\| =O(1/n^{1-\alpha})$.

For the opposite direction, by Lemma \ref{growth-eht} we obtain $\overline{(I-T)X} \subset D(H)$.
Hence $D(H)=\overline{(I-T)X}$, so by Theorem \ref{becker}(ii) $T$ is uniformly ergodic, and 
$(I-T)X$ is closed. By Theorem \ref{C-alpha} $T$ is $(C,\alpha)$ uniformly ergodic.
\end{proof} 
\medskip

Theorem \ref{DL} uses \cite{DL}, which is valid for $X$ over $\mathbb R$ or $\mathbb C$.
When $X$ is over $\mathbb C$, the operators  $(I-T)^\beta$, $0<\beta <1$, were defined in
\cite{AGL} also when $T$ is  only  $(C,\alpha)$-bounded, and it was noted that one gets a
$C_0$-semigroup $\{(I-T)^s\}_{s\ge 0}$ on  $\overline{(I-T)X}$. It was shown in 
\cite[Proposition 4.8]{AGL} that when $T$ is $(C,\alpha)$-bounded, for $0< \beta<1$ we have
$\overline{(I-T)^\beta X}=\overline{(I-T)X}$. In order to extend Theorem \ref{DL} to 
$(C,\alpha)$-bounded operators, we need the following  Proposition.

\begin{prop}\label{different}
Let $0<\alpha<1$ and let $T$ be $(C,\alpha)$-bounded on a complex Banach space $X$. 
Assume that $X=\overline{(I-T)X}$ and $(I-T)X$ is not closed; then 
	
(i) $(I-T)^\beta X$ is not closed for any $0<\beta <1$.

(ii) $(I-T)^\gamma X  \supsetneq (I-T)^\beta X$ for $0<\gamma < \beta  \le 1$.
\end{prop}
\begin{proof}
(ii)  Since $X=\overline{(I-T)X}$, by Lemma \ref{me} $T$ is mean ergodic, with $F(T)=\{0\}$. 
Note that if $(I-T)^\gamma y= 0$ for some  $0<\gamma<1$, then by the semigroup property 
$(I-T)y=(I-T)^{1-\gamma}(I-T)^\gamma y=0$, so $y=0$.

 Let $0<\gamma < \beta  <1$ satisfy $(I-T)^\gamma X = (I-T)^\beta X$. 
Since for any $\delta \in (\gamma,\beta)$ we have 
$(I-T)^\beta X \subset (I-T)^\delta X \subset (I-T)^\gamma X$, we may assume 
$\beta-\gamma <1-\alpha$. Then for $x \in X$
there is $y \in X$ such that $(I-T)^\gamma x=(I-T)^\beta y$, which yields, with 
$\delta=\beta -\gamma$, that $(I-T)^\gamma x=(I-T)^\gamma (I-T)^\delta y$.
 Hence $(I-T)^\gamma[x-(I-T)^\delta y]=0$, which shows $x=(I-T)^\delta y$. We conclude that
$(I-T)^\delta X=X$. By \cite[Theorem 8.4]{AGL} and \cite[Lemma 2.19]{DL} we have 
$$
\overline{(I-T)X}=X=(I-T)^\delta X \subset D(H) \subset \overline{(I-T)X}.
$$
 Now by Theorem \ref{becker}(ii), $T$ is uniformly ergodic, so $(I-T)X$ is closed, a contradiction.
\smallskip

(i) If $(I-T)^\beta X$ is closed, then $(I-T)^\beta X=\overline{(I-T)X}$ by 
\cite[Proposition 4.8]{AGL}. By the semigroup property,  for $0<\gamma <\beta$ we have
$$
\overline{(I-T)X} =(I-T)^\beta X \subset (I-T)^\gamma X \subset \overline{(I-T)X}.
$$
Hence $(I-T)^\gamma X = (I-T)^\beta X$, which contradicts (ii).
\end{proof}

\begin{thm} \label{c-alpha}
Let $0<\alpha<1$ and let $T$ be $(C,\alpha)$-bounded on a complex Banach space $X$. 
Then $T$ is uniformly ergodic if and only if $D(H)=(I-T)X$.
\end{thm}
\begin{proof} Assume $D(H)=(I-T)X$. By Lemma \ref{me}, $\frac1n \|T^n\| \to 0$, so we have 
to prove only that $(I-T)X$ is closed. Denote $Y=\overline{(I-T)X}$, which is 
$T$-invariant, and put $S:=T_{|Y}$. Obviously also $S$ is $(C,\alpha)$-bounded on $Y$.
By Lemma \ref{me}, $x \in Y$ if and only if $\frac1n\sum_{k=1}^n T^k x \to 0$, so we have that
$Y=\overline{(I-S)Y}$. If we show that $(I-S)Y=Y$, then 
$\overline{(I-T)X}=Y=(I-S)Y\subset (I-T)X$, which will yield that $(I-T)X$ is closed.

By \cite[Proposition 4.8]{AGL}, for any $0<\beta<1$ we have 
$$
(I-T)X \subset (I-T)^\beta X \subset \overline{(I-T)^\beta X}=\overline{(I-T)X}=Y.
$$
Hence, for $0< \gamma <\beta <1$ we have, by the semigroup property (see \cite[Section 4]{AGL})
$$
(I-T)^\beta x = (I-T)^\gamma (I-T)^{\beta-\gamma} x \in (I-S)^\gamma Y,
$$
so $(I-T)^\beta X \subset(I-S)^\gamma Y$ for $0<\gamma<\beta<1$.

Assume that $(I-S)Y$ is not closed; fix $0< \beta<1-\alpha$ and $0<\gamma<\beta$.
Applying Proposition \ref{different}(ii) to $S$, we obtain that 
	$(I-S)^{(\gamma+\beta)/2} Y \subset (I-S)^\gamma Y$ with no equality. 
	Let $y \in (I-S)^\gamma Y$ which is not in $(I-S)^{(\gamma+\beta)/2} Y$.  Since 
$\overline{(I-S)Y}=Y$, \cite[Theorem 8.4]{AGL} yields convergence of 
$\sum_{k=1}^\infty \frac1{k^{1-\gamma}}S^k y$, and by \cite[Lemma 2.19]{DL} 
$\ \sum_{k=1}^\infty \frac1{k}S^k y$ converges, so $y \in D(H)$; then by assumption 
$y \in (I-T)X$.
	But then $y \in (I-T)^\beta X \subset (I-S)^{(\beta+\gamma)/2}Y$, a 
contradiction. Hence $(I-S)Y$ is closed, and since $S$ is mean ergodic on $Y$ and has no fixed 
points, $(I-S)Y=Y$. This implies that $(I-T)X$ is closed, hence  $T$ is uniformly ergodic.
\smallskip

The converse is Proposition \ref{coboundaries}.
\end{proof}

\begin{cor}
Let $T$ on a complex Banach space $X$ be $(C,\alpha)$ bounded, $0<\alpha< 1$. 
If $D(H)=(I-T)X$, then $T$ is $(C,\beta)$ uniformly ergodic for any $\beta > \alpha$.
If in addition $\|T^n\|/n^\alpha \to 0$, then $T$ is $(C,\alpha)$ uniformly ergodic.
\end{cor}
\begin{proof}
By Theorem \ref{c-alpha}, $T$ is uniformly ergodic. Hence ${(I-T)X}$ is closed, 
and $X={F(T)\oplus (I-T)X}$.  
The uniform convergence of the $(C,\beta)$ averages  is by \cite[Lemma 1]{De}.

When $\|T^n\|/n^\alpha \to 0$, Theorem \ref{C-alpha} yields $(C,\alpha)$ uniform ergodicity.
\end{proof}

\bigskip

\section{Rotational uniform ergodicity}

A very simple example in \cite{L2} shows that $T$ may be power-bounded and uniformly ergodic while
 $-T$ is not. In that example also $T^2$ is not uniformly ergodic. However, $\|T^nx\| \to 0$ for 
every $x$. A different example was presented by Lyubich and Zem\'anek \cite[Example 3]{LZ},
who gave necessary and sufficient spectral conditions for relative compactness in operator norm of 
$\{M_n(T): n\ge 1\}$.
\smallskip

{\bf Definition.} We call an operator $T$ on a complex Banach space $X$   {\it rotationally uniformly
ergodic} if $\lambda T$ is uniformly ergodic for every $\lambda$ in the unit circle $\mathbb T$.
\smallskip

If $T$ on $X$ complex satisfies $\|T^n\| \to 0$, then it is rotationally uniformly ergodic.

If $T$ is quasi-compact (i.e. $\|T^n-K\|<1$ for some $n \ge 1$ and $K$  compact)
with $n^{-1}T^n \to 0$ in the weak operator topology, then $T$ is uniformly ergodic
\cite[Theorem VIII.8.4]{DS}, and $n^{-1}\|T^n\| \to 0$. Since $\lambda T$ is quasi-compact 
for every $\lambda \in \mathbb T$, $\ T$ is rotationally uniformly ergodic.

\begin{thm} \label{rotational}
The following are equivalent for $T$ on a complex Banach space $X$:

(i) $T$ is rotationally uniformly ergodic.
 
(ii) $r(T) \le 1$ and every $\lambda \in \sigma(T)\cap\mathbb T$ is a simple pole of the resolvent.

 (iii) The set $\{T^n: n\ge 0\}$ is relatively compact in operator norm.
\smallskip

\noindent
If either of the above conditions holds, $T$ is power-bounded and $\sigma(T)\cap \mathbb T$ is finite.  
\end{thm}
\begin{proof}

(i) implies (ii). Let $T$ be rotationally uniformly ergodic.  Clearly $r(T) \le 1$, since 
	$n^{-1}\|T^n\| \to 0$ by uniform ergodicity. If $r(T) <1$ (ii) and (iii) hold trivially.
	Assume $r(T)=1$.  When $\lambda \in \sigma(T)\cap\mathbb T$, the uniform ergodicity of 
	$\bar\lambda T$ yields by \cite[Theorem 3.16]{D} that 1 is a simple pole 
of the resolvent of  $\bar\lambda T$, so $\lambda$ is a simple pole  of the resolvent of $T$.

(ii) implies (iii). If $r(T)<1$, then $\|T^n\| \to 0$ and (iii) obviuosly holds.
Assume now that $r(T)=1$ and  every $\lambda\in\sigma(T)\cap\mathbb T$ is a simple pole of the resolvent.
The spectral points in $\mathbb T$, being poles, are isolated in $\sigma(T)$, so there are necessarily
 only finitely  many of them. Then we obtain (iii) from \cite[Theorem 2]{Sw} or \cite[Theorem 3]{KW} 
(where the next-to-last line of the proof should be replaced by {\it The sequence 
$\{(\lambda_1,\dots,\lambda_k)^n\}_{n\ge 1}$ in  $\mathbb T^k$ is conditionally compact}).

(iii) implies (i). Fix $\lambda \in \mathbb T$. By (iii) also $\{(\lambda T)^n : n\ge 0\}$ is relatively 
compact in norm, and we apply \cite[Theorem 1]{KW} to $\lambda T$ and obtain that
$\lambda T$ is uniformly ergodic.
%
%
\smallskip

Power-boundedness follows from (iii).	Finiteness of $\sigma(T)\cap\mathbb T$ follows from (ii). 
\end{proof}

{\bf Remarks.} 1. For $T$ power-bounded, the equivalence of (ii) and (iii) follows also from Lemma 
2.12 and Proposition 2.13 of \cite{DG}.  Power-boundedness is not assumed in \cite{KW} or \cite{Sw}.

2. A slightly different proof of (iii) implies (i) is by observing (as in \cite{KW}) that by Mazur's 
theorem \cite[Theorem V.2.6]{DS}, (iii) implies that $\{M_n(T): n\ge 1\}$ is also relatively compact.
By \cite[Theorem 4]{LZ}, $1$ is at most a simple pole of the resolvent, and since 
$n^{-1}\|T^n\| \to 0$ by (iii), $T$ is uniformly ergodic by \cite{D}. We then apply this to 
$\lambda T$, for any $\lambda \in \mathbb T$.

3. Note that uniform ergodicity does not imply power-boundedness \cite{LSS1}.

4. If $T$ is quasi-compact with $n^{-1}T^n \to 0$ in the weak operator topology, then $T$ is
 power-bounded.

\begin{cor} \label{powers}
Let $T$ be a rotationally uniformly ergodic operator. Then for every 
$\ell \in\mathbb N$ the operator $T^\ell$ is rotationally uniformly ergodic.
\end{cor}
\begin{proof} By Theorem \ref{rotational}, $\{T^n: n\ge 0\}$ is relatively compact
in operator norm; hence also $\{T^{\ell n}: n\ge 0\}$ is relatively compact, which implies
rotational uniform ergodicity of $T^\ell$.
\end{proof}


\begin{thm} \label{all-circle}
The following are equivalent for $T$ on a complex Banach space $X$:

(i) $T$ is uniformly ergodic with $\sigma(T)\cap \mathbb T = \emptyset$.

(ii)  $r(T) \le 1$ and $\sigma(T)\cap \mathbb T = \emptyset$.

(iii) $\|T^n\| \to 0$.

(iv) $\|T^n\|/n \to 0$, and for every $x \in X$, the rotated one-sided ergodic Hilbert transform series  
\begin{equation} \label{rotated-eht}
\sum_{n=1}^\infty \frac{\lambda^n T^n x}n
\end{equation}
 converges in norm for every $\lambda \in \mathbb T$.

(v) $T$ is rotationally uniformly ergodic with $\sigma(T)\cap \mathbb T = \emptyset$.
\end{thm}
\begin{proof}
Clearly (i) implies  $\|T^n\|/n \to 0$, which yields $r(T) \le 1$,  so (ii) holds. 
Obviously (ii) implies (iii).

The convergence (iii) implies that for some $C>0$ and $\rho \in (0,1)$ we have $\|T^n\| \le C\rho^n$:
Indeed, let $\epsilon <1$ and fix $k$ with $\|T^k\| < \epsilon$. For $n>k$ put $r=[n/k]$; then
$$
 \|T^n\| \le \|T^{rk}\| \le \epsilon^r \le \epsilon^{\frac{n}k-1} =\epsilon^{-1} (\epsilon^{1/k})^n.
$$
This yields that the series $\sum_{n=1}^\infty \frac{\|T^n\|}n$ converges in operator norm,
which yields convergence of \eqref{rotated-eht} for every $x \in X$ and every $\lambda \in \mathbb T$.

Assume (iv). Fix $\lambda \in \mathbb T$. The convergence of \eqref{rotated-eht} for every $x$
implies (by Kronecker's Lemma) that $\lambda T$ is mean ergodic, so 
$\sup_n \|M_n(\lambda T)\| < \infty$. We also have that  $D(H_{\lambda T})=X $ is closed, so
applying Theorem \ref{becker}(iv) to $\lambda T$ we obtain that $\lambda T$ is uniformly ergodic.
Hence $T$ is rotationally uniformly ergodic, so by Theorem \ref{rotational} $\sigma(T)\cap \mathbb T$
contains only eigenvalues.
If $\lambda_0 \in \sigma(T)\cap \mathbb T$ is an eigenvalue with eigenvector $x_0 \ne 0$,
then  \eqref{rotated-eht} with $x=x_0$ and $\lambda=\bar\lambda_0$ does not converge, contradicting 
(iv). Hence $T$ has no unimodular eigenvalues. Consequently  (v) holds. Trivially (v) implies (i).
\end{proof}

\begin{thm} \label{one-point}
The following are equivalent for $T$ on a complex Banach space $X$.

(i) $T^n$ converges uniformly.

(ii) $T$ is uniformly ergodic and $\sigma(T)\cap\mathbb T \subset\{1\}$.

(iii) $T$ is rotationally uniformly ergodic and $\sigma(T)\cap\mathbb T \subset\{1\}$.

(iv)  $n^{-1}\|T^n\| \to 0$, and for each $x \in \overline{(I-T)X}$ the series \eqref{rotated-eht}
 converges for every $\lambda \in\mathbb T$.
\end{thm}
\begin{proof} (i) implies uniform ergodicity, and if $T^n \to E \ne 0$ uniformly, then by the
	ergodic decomposition $\sigma(T)\cap\mathbb T=\{1\}$. Hence (ii) holds.

Assume (ii). If $r(T)<1$, then $\|T^n\| \to 0$ and $T$ is rotationally uniformly ergodic.  
If $\sigma(T)\cap\mathbb T=\{1\}$, then by \cite[Theorem 4]{Lu} we obtain norm convergence of $T^n$,  
so (iii) holds by Theorem \ref{rotational}(iii).

Assume (iii). Denote $Y=\overline{(I-T)X}$. Since $T$ is uniformly ergodic, $D(H_T)=Y$ by Proposition 
\ref{coboundaries}. Thus the series \eqref{rotated-eht} converges for $x \in Y$ and $\lambda=1$. Fix
$1 \ne \lambda_0 \in \mathbb T$. Then $\bar\lambda_0 I-T$ is invertible, and $D(H_{\lambda_0 T})=X$,
so the series \eqref{rotated-eht} converges for  every $x \in X$ and $\lambda=\lambda_0$. Thus (iv) holds.

Assume (iv). Let $Y$ be as above. Then $D(H_T)=Y$, so $T$ is uniformly ergodic by Theorem \ref{becker}(ii),
and $X=F(T)\oplus Y$. Put $S=T_{|Y}$. By (iv) the series \eqref{rotated-eht} converges  with $T$
replaced by $S$, for  every $x \in Y$ and $\lambda \in \mathbb T$. By Theorem \ref{all-circle}(iv), 
$\|S^n\| \to 0$, so $r(S)<1$, and by the ergodic decomposition $T^n$ converges in  norm.
\end{proof}

{\bf Remark.}   See also Proposition \ref{kalton}. Additional conditions equivalent  to the norm 
convergence of $T^n$ were obtained by Jachymski \cite[Theorem 5]{Ja}.

\begin{prop} \label{eht}
Let $T$  on a complex Banach space $X$ satisfy $n^{-1}\|T^n\| \to 0$. Assume that
 $\Lambda:= \sigma(T)\cap \mathbb T$ is non-empty. Then:

(i)  For every $x \in X$ and every $\lambda \in \mathbb T$ with $\bar\lambda \notin \Lambda$, 
 the series \eqref{rotated-eht} converges.

(ii) Assume  $T$ is power-bounded. If $\Lambda =\{\lambda_1,\dots,\lambda_k\}$ is finite,
then for every $x $ in $\cap_{j=1}^k (\lambda_j I-T)X$
and every $\lambda \in \mathbb T$ the series \eqref{rotated-eht} converges. 
\end{prop}
\begin{proof}
(i)  If $|\lambda|=1$ and $\bar\lambda \notin \Lambda$, then $ I-\lambda T$ is invertible, 
so $\lambda T$ is uniformly ergodic, and the series converges for every $x \in X=(I-\lambda T)X$, 
by  applying Theorem  \ref{becker}(ii) to $\lambda T$.

(ii) We have to prove the convergence of \eqref{rotated-eht} also when $\lambda=\bar\lambda_j$.
The assumption on $x$ implies $x \in (\lambda_j I-T)X= (I-\lambda T)X \subset D(H_{\lambda T})$,
 by \cite[Proposition 2.2]{B2} applied to (the power-bounded) $\lambda T$.
\end{proof}

{\bf Remark.} Convergence of the series \eqref{rotated-eht}, when $X$ is a Hilbert space and 
$T$ a contraction, was studied in \cite{CCC}.
\smallskip

{\bf Example.} {\it A power-bounded $T$ with $\sigma(T)=\{1\}$, which is not uniformly ergodic.}

Let $V$ be the Volterra operator on $L_2[0,1]$, and put $T:=I-V$. Then $\sigma(T)=\{1\}$, and
$T$ is power-bounded on $X=L_2[0,1]$ by  \cite{A}. If $T$ were uniformly ergodic, 1 would be an
 eigenvalue; but $T$ has no fixed points, so $T$ is not uniformly ergodic. However, $T$ is mean
 ergodic, and  $X=\overline{(I-T)X}$; in fact, by the Esterle-Katznelson-Tzafriri theorem \cite{KT},
$\|T^n(I-T)\| \to 0$, which yields $\|T^n f\| \to 0$ for every $f \in X$. Note that for
$1\ne \lambda \in \mathbb T$ we have $(I-\lambda T)$ invertible, so $\lambda T$ is uniformly
ergodic. By Proposition \ref{eht}(ii), for $x \in (I-T)X$ the series  \eqref{rotated-eht} converges for 
every $\lambda \in \mathbb T$. Since $T$ is not uniformly ergodic, $D(H) \ne \overline{(I-T)X}=X$
by Theroem \ref{becker}, so for some $x$ the series \eqref{rotated-eht} with $\lambda=1$
does not converge.

\begin{thm}
Let $T$ on a complex Banach space $X$ be  rotationally uniformly ergodic, with 
$\sigma(T)\cap \mathbb T$ non-empty, consisting of the poles $\lambda_1,\dots,\lambda_k$. Then:

(i) For every $x \in X$ and every $\lambda \in \mathbb T$ with $\bar\lambda \notin \Lambda$, 
the series \eqref{rotated-eht} converges.

(ii) For fixed $x \in X$, the series \eqref{rotated-eht} converges for {\em every} 
$\lambda \in \mathbb T$ if and only if $x $ is in $\cap_{j=1}^j (\lambda_j I-T)X$.

(iii)  $[\prod_{j=1}^k (\lambda_jI-T) ]X$ is a closed subspace which equals 
${ \cap_{j=1}^k (\lambda_jI-T) X}$.
\end{thm}
\begin{proof} 
Since $T$ is uniformly ergodic, $n^{-1}\|T^n\| \to 0$, and (i) follows from  Proposition \ref{eht}(i).
\smallskip

 (ii) Let $x \in { \cap_{j=1}^k (\lambda_jI-T) X}$.
 By (i), the series \eqref{rotated-eht} converges for each $\lambda \in \mathbb T$ 
such that $\bar\lambda \ne \lambda_j$, $1\le j \le k$. Fix $j$, and let $\lambda=\bar\lambda_j$. 
Since $x \in {(I-\bar\lambda_jT)X}$ and $\bar\lambda_j T$ is uniformly ergodic,
by Proposition \ref{coboundaries} $x \in D(H_{\bar\lambda_j T})$, which means that the series
\eqref{rotated-eht} converges for $\lambda=\bar\lambda_j$. This part of the proof does not
require power-boundedness, unlike Proposition \ref{eht}(ii).

Now fix $x$  such that the series \eqref{rotated-eht} converges for every $\lambda \in\mathbb T$. 
Fix $j,\quad 1\le j \le k$. Since $\bar\lambda_j T$ is uniformly ergodic, by Proposition 
\ref{coboundaries} $x \in D(H_{\bar\lambda_j T})= (I-\bar\lambda_j T)X=(\lambda_j I-T)X$.
\smallskip

(iii) By Theorem \ref{rotational} $T$ is power-bounded.
	By uniform ergodicity of $\bar\lambda_1 T$, $Y_1:={(\lambda_1 I-T)X}=(I-\bar\lambda_1 T)X$ 
is closed, and obviously $T$-invariant. Then $\bar\lambda_2 T_{|Y_1}$ is uniformly ergodic
on $Y_1$, so $Y_2:=(\lambda_2I-T)Y_1=(\lambda_2I-T)(\lambda_1 I-T)X$ is closed, and so on.

The inclusion $[\prod_{j=1}^k (\lambda_jI-T)]X \subset  \cap_{j=1}^k (\lambda_jI-T) X$  is obvious.

Let $X_j:=\{y: Ty=\lambda_j y\}$ be the eigenspace corresponding to $\lambda_j$,
which is the fixed space of the uniformly ergodic $\bar\lambda_j T$, and denote by 
$E_j:=\lim_n M_n(\bar\lambda_j T)$ the ergodic projection onto $X_j$. Since $T$ is 
power-bounded, we can apply \cite[Theorem 2.4]{CL1} with $T_j=\bar\lambda_j T$, $1\le j \le k$, 
and obtain
\begin{equation} \label{decomposition}
X=\overline{\sum_{1\le j\le k} X_j} \oplus {[\prod_{1\le j \le k}(I-\bar\lambda_j T)] X}.
\end{equation}
It is easily seen that $E_jE_\ell=0$ for $j \ne \ell$, so $E:=\sum_{j=1}^kE_j$ is a projection onto
the left hand side summand in \eqref{decomposition}, which vanishes on the right hand side.
Fix $x \in  \cap_{j=1}^k (\lambda_jI-T) X$. By (ii) $\sum_{n=1}^\infty \bar\lambda_j^nT^n x/n$
converges, so $\sum_{n=1}^\infty \bar\lambda_j^nT^n E_jx/n$ converges, which yields that 
$E_jx =0$. Hence $Ex=0$, so $x \in { [\prod_{j=1}^k (\lambda_jI-T) ]X}$.
\end{proof}

{\bf Remark.}  The  points $\lambda_1,\dots,\lambda_k$  are poles, hence eigenvalues. 
When $\bar\lambda =\lambda_j \in \Lambda$, for $x_j \ne 0$ an eigenvector for 
$\lambda_j$ we have $\lambda^n T^nx_j=x_j$, so the series \eqref{rotated-eht} does not converge
when $x=x_j$ and $\lambda=\bar\lambda_j$.

\begin{thm}
Let $T$ be power-bounded, and assume that $\lambda_1,\dots,\lambda _k \in \mathbb T$ are
 such that each $\bar\lambda_j T$ is mean ergodic (e.g. $X$ is reflexive). 
If for every $x \in X_0:=\overline{ [\prod_{j=1}^k (\lambda_jI-T) ]X}$ the series \eqref{rotated-eht}
converges for every $\lambda \in \mathbb T$, then $T$ is rotationally uniformly ergodic, with
$\sigma(T)\cap \mathbb T \subset\{\lambda_1,\dots,\lambda_k\}$.
\end{thm}
\begin{proof}
Let $T_0$ be the restriction of $T$ to the invariant closed subspace $X_0$.  By assunption, for every
$x \in X_0$ the series \eqref{rotated-eht} converges for every $\lambda \in \mathbb T$. 
Hence, by Theorem \ref{all-circle}, $T_0$ is rotationally uniformly ergodic, with 
$\sigma(T_0)\cap \mathbb T=\emptyset$.

Let $X_j=\{y: Ty=\lambda_j y\}$ be the eigenspace of $\lambda_j$ (maybe $\{0\}$), which is
$F(\bar\lambda_j T)$. Since each $\bar\lambda_j T$ is mean ergodic, we can apply 
\cite[Theorem 2.4]{CL1} and obtain $X=X_0 \oplus \overline{\sum_{j=1}^k X_j}$. The ergodic
projections $E_jx :=\lim M_n(\bar\lambda_j T)x$ project onto $X_j$, and satisfy $E_jE_\ell=0$ for 
$j\ne \ell$. Hence $E:= \sum_{j=1}^k E_j$ is a projection onto $\sum_{j=1}^k X_j$, which
is therefore closed, with $EX_0=\{0\}$. Let $T_j:=T_{|X_j}$.Then $T_j=\lambda_j I_{X_j}$ is
 rotationally uniformy ergodic with $\sigma(T_j)\cap\mathbb T =\{\lambda_j\}$ when $X_j \ne \{0\}$.
Putting everything together, we obtain the theorem.
\end{proof}

{\bf Remarks.} 1. When $k=1$, we obtain the result, assuming only $n^{-1}\|T^n\| \to 0$ and 
without requiring $\bar\lambda_1 T$ to be mean ergodic, by applying Proposition \ref{one-point} 
to  $\bar\lambda_1 T$.

2. Assume $k \ge 2$, and denote $Z:=\overline{\cap_{j=1}^k(\lambda_j I-T)X}$. Obviously, in
the above theorem $X_0 \subset Z$. But for $x \in Z$ the proof shows $Ex=0$, so $x \in X_0$;
thus $Z=X_0$.

\begin{cor}
Let $T$ be power-bounded on a reflexive complex Banach space $X$, with spectral radius $r(T)=1$.  
 Then $T$ is rotationally uniformly ergodic if and only if there exist finitely many points
$\lambda_1,\dots,\lambda _k \in \mathbb T$, such that for every 
$x \in \overline{ [\prod_{j=1}^k (\lambda_jI-T) ]X}$ the series \eqref{rotated-eht} converges 
for every $\lambda \in \mathbb T$.
\end{cor}
{\bf Remarks.} 1. When $r(T)<1$ we have $\|T^n\| \to 0$, and Theorem \ref{all-circle} applies,
even without reflexivity.

2. The assumption of power-boundedness is not restrictive, since by Theorem \ref{rotational}
rotational uniform ergodicity implies power-boundedness.

\medskip

\section{Appendix: uniform ergodicity of positive isometries}

\begin{thm} \label{lattice}
Let $T$ be a positive isometry on an infinite-dimensional complex Banach lattice $L$.
If $F(T)$ is finite-dimensional, then $T$ is not uniformly ergodic.
\end{thm}
\begin{proof} Assume that $T$ is uniformly ergodic. Since $F(T)$ is assumed finite-dimensional,
$T$ is quasi-compact, by \cite{L3}. Bartoszek \cite[Theorem 2]{Ba} proved that if a positive 
contraction $T$ is quasi-compact, then there exist $d\ge 1$ and a finite-dimensional 
projection $E$ such that $\|T^{nd}-E\| \to 0$ as $n \to \infty$. Thus $T^d$ is mean ergodic 
and $E= \lim_n M_n(T^d)$ projects onto $F(T^d)$, which is then finite-dimensional.
It follows from the mean ergodic theorem that $\|T^{nd}f\| \to  0$ for $f\in\overline{(I-T^d)L}$; 
but since $T$ is an isometry, $f=0$. Hence $T^d=I$, which means that $F(T^d)=L$, -- 
a contradiction, since $L$ is infinite-dimensional and $F(T^d)$ is finite-dimensional. 
Thus $T$ is not uniformly ergodic.
\end{proof}

{\bf Remark.}  Theorem \ref{lattice} may fail when $F(T)$ is infinite-dimensional; e.g. $T=I$.

\begin{cor}
Let $T$ be as in Theorem \ref{lattice}. Then 
$\ (I-T)L \underset{\ne} \subset D(H) \underset{\ne} \subset \overline{(I-T)L}$.
\end{cor}
\begin{proof} If we had equality in one of the inclusions, then $T$ would be uniformly ergodic,
by Theorem \ref{DL} or by Theorem \ref{becker}(ii).
\end{proof}

\begin{cor} \label{no-inverse}
Let $\theta $ be an ergodic measure preserving transformation on an infinite-dimensional
probability space $(\mathcal S, \Sigma, \mu)$ and for fixed $1 \le p < \infty$ let 
$Tf:=f \circ \theta$ be the induced isometry on $L_p(\mu)$. Then $T$ is not uniformly ergodic, 
and there exists $f \in \overline{(I-T)L_p}$ such that the series $\sum_{n=1}^\infty \frac1n T^nf$
does not converge in $L_p$ norm.
\end{cor}

{\bf Remarks.} 1. In contrast to Proposition \ref{halmos}, in Corollary  \ref{no-inverse}
invertibility is not assumed.

2. The proof of Proposition \ref{halmos} relies on \cite{Ha}, where the proof uses
approximation of $T$ by a unitary operator induced by a periodic transformation. Corollary
\ref{no-inverse} follows from quasi-compactness arguments in Theorem \ref{lattice}, 
which do not use the one-sided ergodic Hilbert transform.

3. For an ergodic theoretic proof of Corollary \ref{no-inverse}, we note that 
Example 2 of Dowker and Erd\"os \cite{DE} holds also for non-invertible transformations, when in 
its proof we replace Rokhlin's lemma by its non-invertible version, due to Heinemann and Schmitt 
\cite{HS}.
The result of \cite{DE} yields the existence of $f \in L_\infty(\mu)$ with integral zero such that  
the series does not converge in any $L_p$ norm, $1 \le p< \infty$, and we apply 
Proposition \ref{coboundaries} to conclude that $T$ is not uniformly ergodic in $L_p(\mu)$.
\smallskip

The following is a generalization of Proposition \ref{halmos}; ergodicity is replaced
by aperiodicity, and instead of measure preserving, $\theta$ is assumed only {\it non-singular}
-- $\mu\theta^{-1}$ is equivalent to $\mu$.

\begin{prop} \label{nadkarni}
Let $\theta$ be an invertible non-singular bi-measurable aperiodic transformation of an atomless 
probability space $(\mathcal S,\Sigma,\mu)$, and define on $L_2(\mu)$ the operator
$$
	Tf=\sqrt{\frac{d(\mu\theta^{-1})}{d\mu}} \cdot f\circ\theta^{-1}.
$$
Then $T$ is unitary on $L_2(\mu)$ and is not uniformly ergodic.
\end{prop}
\begin{proof} By the change of variables formula, $T$ is unitary. By \cite[Theorem 3.4]{Na},
$\sigma(T)=\mathbb T$. Hence $1 \in \sigma(T)$ is not isolated in the spectrum, 
so is not a pole of the resolvent. By \cite[Theorem 3.16]{D} $T$ is not uniformly ergodic.
\end{proof}

{\bf Remarks.} 1. If we replace aperiodicity by ergodicity, then $F(T)$ is one-dimensional,
and Proposition \ref{nadkarni} follows from Theorem \ref{lattice}.

2.  Corollary \ref{no-inverse} holds also when we replace ergodicity by aperiodicity. 
Adapting the proof of \cite[Theorem 3.4]{Na} and using the non-invertible version of
Rokhlin's lemma \cite{HS}, we can prove, for any $1 \le p< \infty$, that 
$\mathbb T \subset \sigma(T)$ for $T$ acting on $L_p$, so $T$ is not uniformly ergodic.
\medskip

\bigskip

{\bf Acknowledgement.} The authors are grateful to Yuri Tomilov for bringing article \cite{AGL}
to their attention.

\end{document}